\documentclass [a4,10pt,oneside]{article}
\usepackage{amsmath, amssymb, amsthm, amsfonts, amsxtra, latexsym, amscd,
pb-diagram,  graphics,rotating,xy}

\theoremstyle{plain}
\newtheorem{dl}{Theorem}[section]

\newtheorem{dn}[dl]{Definition}

\newtheorem{thm}{\bf Theorem}
\newtheorem{lem}[thm]{\bf Lemma} 
\newtheorem{pro}[thm]{\bf Proposition} 
\newcommand{\A}{\mathcal{A}}
\newcommand{\tx}{\otimes }

\newcommand{\Fm}{\widetilde{F}}

\newcommand{\Gm}{\widetilde{G}}

\newcommand{\C}{\mathcal{C}}

\begin{document}
\setlength{\columnsep}{2cm}\title{ON GR-FUNCTORS BETWEEN GR-CATEGORIES: OBSTRUCTION THEORY FOR GR-FUNCTORS OF THE TYPE $(\varphi,f)$}
\author {Nguyen Tien Quang}

\maketitle
\begin{abstract}
Each Gr-functor of the type $(\varphi,f)$ of a Gr-category of the type $(\Pi,\C)$ has the obstruction be an element $\overline{k}\in H^3(\Pi,\C).$ When this obstruction vanishes, there exists a bijection between congruence classes of Gr-functors of the type $(\varphi,f)$ and the cohomology group $H^2(\Pi,\C).$ Then the relation of Gr-category theory and the group extension problem can be established and used to prove that each Gr-category is Gr-equivalent to a strict one.
\end{abstract}
%\tableofcontents

\section{Introduction and Preliminaries}
\subsection {Introduction}
A monoidal category can be ``refined'' to become a category with group structure when the definition of an {\it invertible object} is added (see [5], [8]). Then, if the background category is a {\it groupoid,} we have the definition of a  monoidal category {\it group-like} (see [2])), or a {\it Gr-category} (see [9]).

The structure of a Gr-category relates closely to the group extension problem and the group cohomology theory. Each Gr-category is Gr-equivalent to its reduced Gr-category. That is a Gr-category of the type $(\Pi,A,\xi),$ where $\Pi$ is a group, $A$ is a $\Pi$-module and $\xi$ is a 3-cocycle of $\Pi,$ with coefficients in $A.$ Therefore, we can classify Gr-categories of the type $(\Pi,A)$ with the cohomology group $H^3(\Pi,A).$ This result has been applied and extended in some studies on graded extension, where Gr-categories are called {\it categorical groups} (see [1]).

In this paper, let us consider the Gr-functor classification problem. This problem can be reduced to consider Gr-functors of Gr-categories of the type $(\Pi,A).$ Each such Gr-functor is a functor of the type $(\varphi,f).$ Conversely, each functor of the type $(\varphi,f)$ induces an element in the cohomology group $H^3(\Pi,A),$ called the {\it obstruction} of the functor of the type $(\varphi,f).$ When this obstruction vanishes, Gr-functors exist and can be classificated with the group $H^2(\Pi,A).$

The definition of the obstruction of Gr-functors relates to the definition of the group extension problem. Applying these 2 theories, we can prove that an arbitrary Gr-category is equivalent to a strict one. (The analogous problem for a {\it monoidal category} can be solved in many different way (see [4])).

\subsection {Elementary concepts}

\noindent Let us start with some elementary concepts of a monoidal category.\\
\indent \emph{ A monoidal category} $(\mathcal{C} ,\otimes, I, a, l,
r)$ is a category $\mathcal{C}$ together with a tensor
product $\otimes: \mathcal{C} \times \mathcal{C}\rightarrow
\mathcal{C};$ an object $I,$ called \emph{the unitivity object} of the
category; and the natural isomorphisms
\begin{gather*}
a_{A, B, C}:A\otimes(B\otimes C)\rightarrow (A\otimes B)\otimes C\\
l_A:I\otimes A\rightarrow A\\
r_A:A\otimes I\rightarrow A \end{gather*} respectively,
called \emph{the associativity constraint, the left unitivity
constraint} and \emph{the right unitivity constraint}. These
constraints must satisfy the Pentagon Axiom
\begin{equation}(a_{A,B,C}\otimes id_D)\  a_{A,B\otimes C,D}\ (id_A\otimes
a_{B,C,D})=a_{A\otimes B,C,D}\ a_{A,B,C\otimes D},\tag{1.1}
\end{equation} and the Triangle Axiom
\begin{equation}
id_A\otimes l_B=(r_A\otimes id_B)a_{A,I,B}.\tag{1.2} \end{equation}
\indent A monoidal category is \emph{strict} if the associativity $a$ and unitivity constraints $\ l,\ r$ are all identities.\\
\indent Let $\mathcal{C}=(\mathcal{C}, \otimes, I, a, l, r)$ and
$\mathcal{C'}=(\mathcal{C'}, \otimes, I', a', l', r')$ be monoidal
categories. \emph{A monoidal functor} from $\mathcal{C}$ to
$\mathcal{C'}$ is a triple $(F, \widetilde{F},\hat{F})$ where
$F:\mathcal{C} \rightarrow \mathcal{C'}$ is a functor, $\hat{F}$ is
an isomorphism from $I'$ to $FI$, and $\widetilde{F}$ is a natural
isomorphism
\begin{equation}
\widetilde{F}_{A,B}:FA\otimes FB\rightarrow F(A\otimes B) \tag{1.3}
\end{equation} satisfying following commutative diagrams
\begin{equation}
\begin{CD}
FA\otimes(FB\otimes FC )@>id\otimes \widetilde{F}>> FA\otimes F(B\otimes C)@>\widetilde{F} >> F(A\otimes (B\otimes C))\\
@ V a' VV              @.           @VV F(a) V \\
(FA\otimes FB )\otimes FC @>\widetilde{F} \otimes id>> F(A\otimes
B)\otimes FC @>\widetilde{F}
>> F((A\otimes B)\otimes C)
\end{CD}
\tag{1.4}
\end{equation}
\begin{equation}
\begin{CD}
FA \otimes FI @ > \widetilde{F} >> F(I \otimes A)\\
@ A id \otimes \hat{F} AA @ VV F(r) V \\
FA \otimes I' @ > r' >> FA
\end{CD}
\tag{1.5a}
\end{equation}
\begin{equation}
\begin{CD}
FI \otimes FA @ > \widetilde{F} >> F(A \otimes I)\\
@ A \hat{F} \otimes id AA @ VV F(l) V \\
I' \otimes FA @ > l' >> FA
\end{CD}
\tag{1.5b}
\end{equation}
\indent \emph{A natural monoidal transformation}
$\alpha:(F,\widetilde{F}, \hat{F})\to(G,\widetilde{G},\hat{G})$
between monoidal functors from $\mathcal{C}$ to $\mathcal{C'}$ is a
natural transformation $\alpha: F\rightarrow G,$ such that the
following diagrams
\begin{equation}
\begin{CD}
FX \otimes FY @ > \alpha_X \otimes \alpha_Y >> GX \otimes GY \\
@ V \widetilde{F} VV @ VV \widetilde{G} V \\
F(X\otimes Y) @ >\alpha_{X\otimes Y} >> G(X\otimes Y)
\end{CD}
\tag{1.6}
\end{equation}
\[
\begin{diagram}
\node[1]{FI}
 \arrow[2]{e,t}{\overline{\alpha}_I}
\node[2]{GI}
\\
\node[2]{I'}\arrow[1]{nw,r,2}{\widehat{F}}
\arrow[1]{ne,r,2}{\widehat{G}}
\end{diagram}
\tag{1.7}
\]
commute, for all pairs $(X,Y)$ of objects in $\mathcal{C}.$\\
\indent \emph{An monoidal isomorphism} is a monoidal transformation as well as a natural isomorphism.\\
\indent \emph{ A monoidal equivalence} between monoidal categories is a monoidal functor $F:\mathcal{C} \to \mathcal{C'},$ such that there exists a monoidal functor $G:\mathcal{C'}\to \mathcal{C}$ and natural monoidal isomorphism $\alpha: G.F\to id_{\mathcal{C}}$ and $ \beta:F.G\rightarrow id_{\mathcal{C'}}$.\\
\indent $\mathcal{C}$ and $\mathcal{C'}$ are \emph{monoidal
equivalent} if there exists a monoidal equivalence between them.
%%%%%%%%%%%%%%%%%%%%%%%%%%%%%%

\subsection {A reduced Gr-category and canonical equivalences}

\indent Let $\mathcal{C}$ be a Gr-category. Then, the set of congruence classes of invertible objects $\Pi_0(\ C)$ of $\C$ is a group with the operations induced by the tensor product in $\C,$ and $\Pi_1(\mathcal{C})$ of automorphisms of the unitivity object $I$ is an abel group with the operation, denoted by +, induced by the composition of the arrows. Moreover, $\Pi_1(\mathcal{C})$ is a $\Pi_0(\mathcal{C})$-module with the
action
$$su=\gamma_{A}^{-1}\delta_{A}(s),\ A\in s,\ s\in \Pi_0(\mathcal{C}),\ u\in \Pi_1(\mathcal{C})$$
where $\delta_A,\ \gamma_A$ are defined by the following commutative
diagrams
$$\begin{CD}
X @>\gamma_X(u) >> X \\
@A l_X AA     @AA l_X A \\
I \otimes X @>u \otimes id >> I \otimes X
\end{CD}
\qquad \qquad\qquad
\begin{CD}
X @>\delta_X(u) >> X \\
@A r_X AA     @AA r_X A \\
X \otimes I @>id \otimes u >> X \otimes I
\end{CD}$$

\indent The reduced Gr-category $\mathcal S$ of a Gr-category
$\mathcal{C}$ is the category whose objects are the elements of
$\Pi_0(\mathcal{C})$ and whose arrows are automorphisms of the
form $(s,u),$ where $s \in \Pi_0(\mathcal{C}),\ u \in
\Pi_1(\mathcal{C})$. The composition of two arrows is induced by
the addition in $\Pi_1(\mathcal{C})$ as follows
$$(s,u).(s,v)=(s,u+v).$$
The category $\mathcal{S},$ equivalent to $\C$ thanks to the canonical equivalence, is built as follows. For each $s=\overline X \in \Pi_0(\mathcal{C}),$ we choose a
representative $X_s\in \mathcal{C}$ and for each $X \in s,$ we
choose an isomorphism $i_X:\ X_s\rightarrow X$ such that
$i_{X_s}=id$. The family $(X_s, i_X)$ is called a \emph{ stick} of the
Gr-category $\mathcal{C}.$ For any $(X_s, i_X),$ we obtain two
functors
\[\begin{cases}
G:\mathcal{C}\rightarrow \mathcal S\\
G(X)=\overline X=s\\
G(X \stackrel
{f}{\rightarrow}Y)=(s,\gamma_{X_s}^{-1}(i_{Y}^{-1}fi_X))
\end{cases}\qquad\qquad
\begin{cases}
H:\mathcal S\rightarrow \mathcal{C}\\
H(s)=X_s\\
H(s,u)=\gamma_{X_s}(u)
\end{cases}\]
\indent The functors $G$ and $H$ are equivalences of categories together with
the natural transformations
$$\alpha=(i_{X}): HG\cong id_{\C}\ ;\ \beta=id:GH\cong id_\mathcal S$$
They are called \emph{ canonical equivalences}.\\
\indent With the structure conversion (see [8]) by the quadruple $(G, H, \alpha,
\beta),$ $\mathcal S$ can be equiped with the following operation to become a Gr-category
\begin{gather*}
s\otimes t=s.t,\quad s,t\in\Pi_0(\mathcal{C})\\
(s,u)\otimes(t,v)=(st,u+sv),\quad u,v\in\Pi_1(\mathcal{C}).
\end{gather*}
The Gr-category $\mathcal S$ has the unitivity constraint be strict and the associativity constraint be a normalized 3-cocycle $\xi$ in $\Pi_0(\C)$ with coefficients in $\Pi_0(\C)$-module $\Pi_1(\C).$\\
Moreover, the equivalence $G,\ H$ become Gr-equivalences together with natural transformations
$$\Breve{G}_{A,B}=G(i_A\otimes i_B)\ , \
 \Breve{H}_{s,t}=i_{A_s\otimes A_t}^{-1}$$
The Gr-category $\mathcal S$ is called a \emph{ reduction} of the
Gr-category $\C.$ We also say that $\mathcal S$ is of the type $(\Pi, A, \xi)$ or
simply the type $(\Pi,A)$ when $\Pi_0(\mathcal{C}),
\Pi_1(\mathcal{C})$ are, respectively, replaced with the group $\Pi$ and the
$\Pi$-module $A.$
%%%%%%%%%%%%%%%%%%%%%%%%%%%%%%%%%%%%%%%
\section{
The obtruction and classification of Gr-functors of the type $(\varphi,f)$}
In this section, we will show that each Gr-functor $(F, \Fm):\C\to\C^{'}$ induces a Gr-functor $\overline{F}$ on their reduced Gr-categories, and this correspondence is 1-1. This allows us to study the Gr-functor existence problem and classify them on Gr-categories of the type $(\Pi, A).$ First, we have
\begin{pro} $[9]$
Let $(F, \widetilde{F}): \mathcal{C}\to\mathcal{C'}$ be a Gr-functor.
Then, $(F,\widetilde{F})$ induces the pair of group homomorphisms
\begin{gather*}
 F_0: \Pi_{0}(\C)\rightarrow \Pi_{0}(\C^{'})\ \ ; \ \ \overline A\mapsto \overline {FA}\\
F_1:\Pi_{1}(\C)\rightarrow \Pi_{1}(\C^{'})\ \ ; \ \ u\mapsto \gamma_{FI}(Fu)
\end{gather*}satisfying $F_1(su)= F_0(s)F_1(u).$
\end{pro}
The pair $(F_0,F_1)$ is called \emph{ the pair of induced homomorphisms} of
the Gr-functor $(F,\widetilde{F}).$ Let $\mathcal S, \mathcal S'$
be, respectively, the reduced Gr-categories of $\mathcal{C}, \mathcal{C}'.$ Then, the functor $\overline{F}:
\mathcal S\to \mathcal S'$ given by
$$\overline F (s)=F_0(s), \overline F (s,u)=(F_0 s,F_1 u)$$
is called the {\it induced functor} of $(F,\widetilde{F})$ on reduced Gr-categories.\\
\begin{pro}
Let $\overline F$ be the induced functor of the Gr-functor $(F,
\widetilde{F}): \mathcal{C}\to\mathcal{C}'$. Then the diagram
\[\begin{diagram}
\node{\mathcal{C}}\arrow{e,t}{F}
\node{\mathcal{C}'}\arrow{s,r}{G'}\\
\node{\mathcal S}\arrow{n,l}{H}\arrow{e,t}{\overline F}
\node{\mathcal S'}
\end{diagram}\]
commutes, where $H, G'$ are canonical equivalences, and therefore $\overline{F}$ induces a Gr-functor.
\end{pro}
\indent In order to prove this proposition, we consider the two following
lemmas
\begin{lem} Let $\mathcal{C}, \mathcal{C'}$ be $\otimes$-categories with, respectively, constraints $(I,l,r)$ and $(I',l',r').$ Let $(F,\widetilde{F}):\mathcal{C} \to \mathcal{C'}$ be a $\otimes$-functor which is compatible with the unitivity constraints. Then, the following diagram
\[\begin{diagram}
\node{FI}\arrow{e,t}{\gamma_{_{FI}}(u)}
\node{FI}\\
\node{I'}\arrow{n,l}{\hat{F}}\arrow{e,t}{u}
\node{I'}\arrow{n,r}{\hat{F}}
\end{diagram}\]
commutes, where $\hat{F}$ is the isomorphism induced by $(F,
\widetilde{F})$.
\end{lem}
%\vspace{0.2cm}
\begin{proof} It is clear that $\gamma_{I'}(u)=u$. Moreover, the family $(\gamma_{A'}(u)), \ A' \in \mathcal{C'},$ is a morphism of the identity functor $id_{\mathcal{C'}}.$ So the above diagram commutes.
\end{proof}
\begin{lem} If the assumption of Lemma 3 holds, we have
$$F\gamma_{_A}(u) = \gamma_{_{FA}}(\gamma^{-1}_{_{FI}}Fu).$$
\end{lem}
\begin{proof}
Consider the following diagram
\begin{center}
\setlength{\unitlength}{1cm}
\begin{picture}(12,5)

\put(1.8,1.5){$I'\otimes FA$} \put(4.1,1.5){$FI\otimes FA$}
\put(6.6,1.5){$F(I\otimes A)$} \put(9.6,1.5){$FA$}
\put(1.8,3.3){$I'\otimes FA$} \put(4.1,3.3){$FI\otimes FA$}
\put(6.6,3.3){$F(I\otimes A)$} \put(9.6,3.3){$FA$} {\scriptsize
\put(3.1,1.2){$\hat{F}\otimes id$} \put(8.5,1.2){$F(l_A)$}
\put(6,1.7){$\widetilde{F}$} \put(0.7,2.5){$\gamma^{-1}_{FI}Fu
\otimes id$} \put(4.2,2.5){$Fu\otimes id$} \put(6.7,2.5){$F(u\otimes
id)$} \put(9.9,2.5){$F\gamma_A(u)$} \put(3.2,3.5){$\hat{F}\otimes
id$} \put(8.5,3.5){$F(l_A)$} \put(6,3){$\widetilde{F}$}
\put(6,0.3){$l$} \put(6,4.6){$l$}

\put(3.3,1.6){\vector(1,0){0.8}} \put(5.8,1.6){\vector(1,0){0.8}}
\put(8.3,1.6){\vector(1,0){1.3}} \put(3.3,3.4){\vector(1,0){0.9}}
\put(5.8,3.4){\vector(1,0){0.8}} \put(8.3,3.4){\vector(1,0){1.3}}
\put(2.4,3.1){\vector(0,-1){1.2}} \put(4.8,2.4){\vector(0,-1){0.5}}
\put(7.3,2.4){\vector(0,-1){0.5}} \put(9.7,4.4){\vector(0,-1){0.6}}
\put(9.7,0.7){\vector(0,1){0.6}} \put(9.7,3.1){\vector(0,-1){1.2}}

\put(2.4,1.3){\line(0,-1){0.6}} \put(2.4,0.7){\line(1,0){7.3}}
\put(2.4,3.6){\line(0,1){0.8}} \put(2.4,4.4){\line(1,0){7.3}}
\put(4.8,2.8){\line(0,1){0.3}} \put(7.3,2.8){\line(0,1){0.3}} }
\put(3.2,2.5){(1)} \put(6,2.5){(2)} \put(8.5,2.5){(3)}
\put(6,1){(4)} \put(6,4){(5)}

\end{picture}
\end{center}
In this diagram, the regions (4) and (5) commute thanks to the compatibility
of the functor $(F,\widetilde{F})$ with the associativity constraints. The region (3) commutes thanks to the definition of $\gamma_A$
(image through $F$), the region (1) commutes by Lemma 3. The region (2) commutes thanks to the naturality of the isomorphism
$\widetilde{F}.$ Therefore, the outside commutes, i.e.,
$$F\gamma_A(u)=\gamma_{FA}\big(\gamma^{-1}_{FI}Fu\big).$$
\end{proof}
{\bf The proof of Proposition 2.}
Set $K = G'FH.$ It is easy to verify that $K(s)= \overline{F}(s)$, for $s\in \Pi_0(\C)$. We now prove that $K(s,u)= \overline{F}(s,u)$ for each arrow $u: I \to I:$ 
$$K(u) = G'FH(u) = G'(F \gamma_{A_s}(u)).$$
 \indent Since $H'G' \simeq id_{\mathcal{C'}}$ thanks to natural equivalence $\beta = (i'_{A'}),$ the following diagram
$$\begin{CD}
A_s'  @>i'>>  FA_s \\
@V H'G'F \gamma_{A_s}(u)VV   @VV F \gamma_A(u) V \\
A_s' @>i'>> FA_s
\end{CD}$$
commutes. (Note that $A_s' = H' G'FA_s$).\\
\indent According to Lemma 4, we have
$$ F\gamma_{A_s}(u)=\gamma_{FA_s}(\gamma^{-1}_{FI}Fu).$$
\indent Besides, since the family $(\gamma_{A'})$ is a natural equivalence of the identity functor $id_{\mathcal{C}'}$, the
following diagram
$$\begin{CD}
A_s'  @>i'>>  FA_s \\
@V \gamma_{A_s'}(\gamma^{-1}_{FI}Fu)VV   @VV \gamma_{FA_s}(\gamma^{-1}_{FI}Fu) V \\
A_s'  @>i'>>  FA_s
\end{CD}$$
commutes.\\
Hence $H'G'F \gamma_{A_s}(u)=
\gamma_{A'_s}\big(\gamma^{-1}_{FI}Fu\big)$. From the definition
of $H',$ we have
$$G'F\gamma_{A_s}(u)=\gamma^{-1}_{FI}Fu=F_\circ(u).$$
This means $G'FH=\overline{F}$.\\

\indent Now we describe the Gr-functors on Gr-categories of the type $(\Pi,\A)$ 
\begin{dn}
Let $\mathcal S=(\Pi, A,\xi),\ \mathcal S'=(\Pi',A',\xi')$ be Gr-categories. A functor $F: \mathcal S\to \mathcal S'$ is called a functor of the type $(\varphi,f)$ if
$$F(x)=\varphi(x)\ \ ,\ \ F(x,u)=(\varphi(x), f(u)) $$
and $\varphi:\Pi\to \Pi'$, $f:A\to A'$ is a pair of group homomorphisms satisfying $f(xa)=\varphi (x)f(a)$ for $x\in \Pi, a\in A.$ 
\end{dn}
\begin{thm}
Let $\mathcal S=(\Pi, A, \xi),\ \mathcal S'=(\Pi',A',\xi')$ be
Gr-categories and $(F,\widetilde{F})$ be a Gr-functor from $\mathcal
S$ to $\mathcal S'$. Then, $(F,\Fm)$ is a functor of the type $(\varphi,f).$
\end{thm}
\begin{proof}
\indent For $x,y\in \Pi$, $\widetilde{F}_{x,y}:Fx\otimes Fy\to F(x\otimes y)$ is an arrow in $\mathcal S'$. It follows that $Fx.\ Fy=F(xy),$ so the function $\varphi:\Pi\to\Pi'$, defined by $\varphi(x)=F(x)$ is a group homomorphism on objects.\\
\indent Assume that  $F(x,a)=(Fx,f_x(a))$. Since
$F$ is a functor, we have
$$F((x,a).(x,b))=F(x,a).F(x,b).$$
It follows that
\begin{equation}
 f_x(a+b)=f_x(a)+f_x(b).\tag {2.1}
\end{equation}
So $f_x:A\to A'$ is a group homomorphism for  each $x\in \Pi$. Besides, since $(F,\widetilde{F})$ is a $\otimes$-functor,
the diagram
$$\begin{CD}
Fx.Fy @>\widetilde{F} >> F(xy)\\
@V Fu\otimes Fv VV   @ VV F(u\otimes v) V \\
Fx.Fy @ >\widetilde{F} >> F(xy)
\end{CD}$$
commutes for all $u=(x,a),\ v=(y,b)$. Hence, we have
%\vspace{0.3cm}\\
\begin{gather*}
F(u\otimes v)=Fu\otimes Fv\\
\Leftrightarrow   f_{xy}(a+xb)=f_x(a)+Fx.f_y(b)\\
\Leftrightarrow\  f_{xy}(a)+f_{xy}(xb)=f_x(a)+Fx.f_y(b) \tag {2.2}
\end{gather*}
Applying the relation (2.1) for $x=1,$ we obtain $f_y(a)=f_1(a)$. Thus,
$f_y=f_1$ for all $y\in \Pi.$ Set $f_y=f$ and use (2.1'), we obtain
\begin{equation}
f(xb)=F(x)f(b)=\varphi(x)f(b)\tag{2.3}
\end{equation}
Note that if we regard $\Pi'$-module $A'$ as a $\Pi$-module by the action
$xa'=F(x).a',$ then from (2.1), (2.3), it results $f:A\to A'$ is a
homomorphism of $\Pi$-modules.
\end{proof}
To find the sufficient condition to make a functor of the type $(\varphi,f)$ become a Gr-functor, let us present the definition of {\it obstruction} as in the case of Ann-functors (see [7]).
\begin{dn}
If $F:(\Pi, A,\xi)\to (\Pi', A',\xi')$ is a functor of the type $(\varphi,f),$ then $F$ induces 3-cocycles $\xi_*=f_*\xi,\xi'^*=\varphi^*\xi',$ in which
\begin{gather*}
(f_*\xi)(x,y,z)=f(\xi(x,y,z))\\
(\varphi^*\xi')(x,y,z)=\xi'(\varphi x, \varphi y,\varphi z)
\end{gather*}
The function $k=\varphi^{\ast}\xi'-f_{\ast}\xi$ is called an obstruction of the functor of the type $(\varphi,f).$
\end{dn}
\begin{thm}
The functor $F: \mathcal (\Pi, A,\xi)\to \mathcal
(\Pi',A',\xi')$ of the type $(\varphi,f)$ is a Gr-functor iff the obstruction $\overline{k}=0$ in  $H^3(\Pi, A').$
\end{thm}
\begin {proof}
 If $(F,\widetilde{F}): (\Pi, A,\xi)\to (\Pi',A',\xi')$ is a
Gr-functor, the diagram
$$\begin{CD}
Fx\tx(Fy\tx Fz )@>id \tx \Fm >> Fx\tx F(y\tx z) @> \Fm >>F(x\tx (y\tx z)) \\
@V (\bullet, \xi'_{Fx, Fy, Fz})VV                       @.                                @VV F(\bullet, \xi_{x,y,z})V \\
(Fx\tx Fy)\tx Fz @>\Fm \tx id >>F(x\tx y)\tx Fz @> \Fm >>F((x\tx y)\tx z)
\end{CD}$$
commutes. It implies
$$(Fx, 0)\otimes \Fm_{y,z}+\Fm_{x,yz}+(\bullet, \xi_*(x,y,z))=(\bullet, {\xi'}^*(x,y,z)+\Fm_{x,y}\otimes (Fz, 0)+\Fm_{xy,z}.$$
%%%%%%%%%%%%%%%%%%%%%%%%%%%%%%%%%%%%
%%%%%%%%%%%%%%%%%%%%%%%%%%%%%%%%%%
Since ($\widetilde{F}_{x,y}=(F(xy),k(x,y)):Fx.Fy\to F(xy)$,
where $g:\Pi^2\to A'$ is a function, we may identify $\widetilde{F}$ with $k$ and call $k$ the {\it associated function} with $\Fm.$ Then, we have:
$${\xi'}^*-\xi_*=\delta k,$$
is a 3-coboundary of the group $\Pi$ with coefficients in the $\Pi$-module $A'.$ Therefore, $\overline{\xi'}^*-\overline {\xi_*}=0$  in $H^3(\Pi, A')$.\\
Conversely, from $\overline{\xi'}^*-\overline
{\xi_*}=0,$ there exists a 2-cochain $k\in Z^2(\Pi,A')$ such that
${\xi'}^*-\xi_*=\delta k.$ Take $\Fm_{x,y}$ be associated with $k(x, y),$  we can see that $(F,\widetilde{F})$ is a Gr-functor. 
\end{proof}

Consider the case in which the obstruction of the functor $F$ of the type $(\varphi,f)$ vanishes, we have
\begin{thm}
1) There exists a bijection between the set of congruence classes of Gr-functors of the type $(\varphi,f)$ and the cohomology group $H^2(\Pi, A')$ of the group $\Pi$ with coefficients in the $\Pi$-module $A',$ where $A'$ is a $\Pi$-module with the action $xa'=\varphi(x)a'$.\\
2) If $F:(\Pi,A,\xi)\to(\Pi',A',\xi')$ is a Gr-functor,
there exists a bijection
$$Aut(F)\to Z^1(\Pi, A').$$
\end{thm}
\begin {proof}
1) Suppose that $(F,\widetilde{F}):(\Pi, A, \xi)\to (\Pi', A',
\xi')$ is a Gr-functor, then from Theorem 5, $F=(\varphi, f)$ and 
$${\xi'}^*-\xi_*=\delta k$$
where $k$ is the associated funtion with $\Fm$. Let $k$ be fixed. Now if 
$$(G,\widetilde{G}):(\Pi, A, \xi)\to (\Pi', A', \xi')$$
is a Gr-functor of the type $(\varphi,f),$ then $\Gm\equiv g$ and $${\xi'}^*-\xi_*=\delta g.$$
So $k-g$ is a 2-cocycle. Consider the correspondence:
\[\Phi: class(G,\Gm)\mapsto class(k-g)\]
between the set of congruence classes of Gr-functors of the type $(\varphi,f)$ from $(\Pi, A, \xi)$ to $(\Pi', A', \xi')$ and the group $H^2(\Pi, A').$

First, we will show that the above correspondence is a function. Indeed, suppose that
$$(G',\Gm^{'}):(\Pi, A, \xi)\to (\Pi', A', \xi')$$
\noindent is a Gr-functor of the type $(\varphi,f),$ and $\alpha: G\to G^{'}$ is a Gr-natural transformation, then for all $x,y\in \Pi,$ the diagram
$$\begin{CD}
Gx.Gy @ > \Gm >> G(xy) \\
@ V (\bullet, \alpha_{x})\tx(\bullet, \alpha_y)VV  @ VV(\bullet, \alpha_{xy}) V \\
G^{'}x.G^{'}y  @ > \Gm^{'} >> G^{'}(xy)\\
\end{CD}
$$
commutes. From the definition of the tensor product in the category
$(\Pi', A'),$ we have
$$\alpha_x\otimes \alpha_y=\alpha_x+Gx.\alpha_y.$$
It follows that
$$x\alpha_y-\alpha_{xy}+\alpha_x=g-g',$$
where $g, g'$ are, respectively, associated functions with $\Gm, \Gm'.$

(Since $\widetilde{G}=g,\ \widetilde{G'}=g'$ are 2-cocyles and
$\alpha$ is an 1-cochain, it implies)
\begin{equation}
g-g'=\delta \alpha. \tag{2.4}
\end{equation}
So $\overline {k-g}=\overline {k-g'}\in H^2(\Pi, A').$
We now prove that $\Phi$ is an injection. Suppose that 
$$(G,\Gm),(G',\Gm'):(\Pi, A, \xi)\to (\Pi', A', \xi')$$
\noindent are Gr-functors of the type $(\varphi,f)$ and satisfying
\[\overline {k-g}=\overline {k-g'}\in H^2(\Pi, A').\]
Then, there exists an 1-cochain and such that
\[k-g=k-g'+\delta \alpha\]
i.e. $g'=g+\delta \alpha.$  So the above-mentioned diagram commutes, i.e., $\alpha: G\to G'$ is an $\tx$-morphism. So
\[class (G, \Gm)=class(G', \Gm').\]
Finally, we will show that the correspondence $\Phi$ is a surjection. Indeed, assume that $g$ is an arbitrary 2-cocycle. We have
\[\delta (k-g)=\delta k-\delta g=\delta k = \xi'^*-\xi_*.\]
Then, from Theorem 6, there exist the Gr-functor 
\[(G,\Gm):(\Pi, A, \xi)\to (\Pi', A', \xi')\]
of the type $(\varphi,f)$ and the functorial isomorphism $\Gm=(\bullet,k-g).$ Clearly, $\Phi(G)=\overline{g}.$ So $\Phi$ is a surjection.

2) Let $F=(F,\widetilde{F}):(\Pi, A, \xi)\to (\Pi', A', \xi')$ be a
Gr-functor and $\alpha \in Aut(F).$ Then, the equality (2.4) implies that
$\delta \alpha=0,$ i.e., $\alpha\in Z^1(\Pi, A').$
\end{proof}

\section {The relation between a Gr-category and the group extension problem}
In this section, we built a Gr-category of an abstract kernel and apply it to create the constraints of a Gr-category strict.
\subsection {The Gr-category $\mathcal A_G$ of the group $G$}
For the
given group $G,$ we construct a Gr-category $\mathcal A_G$ whose objects are
elements of the group of automorphisms $Aut(G)$. For the objects $\alpha, \beta$ of $\mathcal A_G,$ we denote
$$ Hom(\alpha, \beta)=\{c\in G|\alpha=\mu_c\circ \beta\}.$$
where $\mu_{c}$ is the automorphism induced by $c.$
For the arrows $c:\alpha\to \beta;\ d:\beta\to\gamma$ of $\mathcal A_G$, their
composition is defined by $d\circ c=d+c$ (the addition in $G$).
Then, the category $\mathcal A_G$ is a strict Gr-category with the tensor
product defined by
\begin{gather*}
\alpha\otimes \beta =\alpha\circ \beta \\
c\otimes d=c+\alpha'(d)
\end{gather*}
where $c:\alpha\to\alpha',\ d:\beta\to\beta'$.\\
The following proposition describes the reduced Gr-category of the Gr-category of an abstract kernel.
\begin{thm}
Let $(\Pi,G,\psi)$ be abstract kernel and $\overline{\xi}\in
H^3(\Pi,Z(G))$ be its obstruction. Let $\mathcal S'=(\Pi',C,\xi')$
be the reduced Gr-category of the strict one $\mathcal A_G.$ Then
$$\Pi'=AutG/IntG\  , \ C=Z(G),$$
and $\psi^*\xi'$ belongs to the cohomology class of $\xi.$
\end{thm}
\begin{proof}
From the definitions of the category $\mathcal A_G,$ and reduced Gr-categories, we have
$$\Pi'=AutG/IntG\ ,\ C=Z(G)$$
We just have to prove that $\overline{\psi^*\zeta'}=\overline{\zeta}.$\\
\indent Indeed, let $(H, \widetilde{H})$ be a canonical
Gr-equivalence from $\mathcal S'$ to $A_G.$ Then, the diagram
\begin{equation}
\begin{CD}
H_r(H_s H_t )@>id\otimes\widetilde{H}>> H_rH_{st}@> \widetilde{H}>>H_{r(st)} \\
@|                       @.                        @VV H(\bullet,\xi'_{r,s,t})V \\
(H_r H_s) H_t @>\widetilde{H}_{r,s}\otimes id >> H_{rs} H_{t} @>
\widetilde{H}>> H_{(rs)t}
\end{CD}
\tag{3.1}
\end{equation}
commutes for all $r, s, t\in \Pi'.$ Since $Au_G$ is strict, we have
$$\gamma_r(u)=u, \ \forall r \in Au_G,\ \forall u \in Z(G)=C.$$
Together with the definition of $H,$ we obtain $H(\bullet,c)=c, \
\forall c \in C.$ The fact that the diagram (3.1) commutes gives us
\begin{equation}
H_r(h_{s,t})+h_{r,st}=-\zeta'_{r,s,t}+h_{r,s}+h_{rs,t}\tag{3.2}
\end{equation}
where $h_{s,t}=\widetilde H_{s,t}$ is a function from $\Pi'\times\Pi' $ to $G.$\\
\indent For the abstract kernel $(\Pi,G,\psi),$ we choose the
function $\varphi=H.\psi:\Pi\to Aut(G).$ Clearly, $\varphi(1)=id_G.$
Moreover, since
$$\widetilde H_{\psi(x), \psi(y)}:H\psi(x)H\psi(y)\to H\psi(xy)$$
is a morphism in $Au_G,$ for all $x,y\in\Pi,$ we have
$$\varphi(x)\varphi(y)=H\psi(x)H\psi(y)=\mu_{f(x,y)}H\psi(xy)=\mu_{f(x,y)}\varphi(xy)$$
where $f(x,y)=\widetilde H_{\psi(x),\psi(y)}.$ So the pair
$(\varphi,f)$ is a factor set of the abstract kernel $(\Pi,G,\psi).$
So, there exists an obstruction $k(x,y,z)\in C=Z(G)$ satisfying
$$\varphi(x)[f(y,z)]+f(x,yz)=k(x,y,z)+f(x,y)+f(xy,z).$$
From the supposition, $\overline{k}=\overline{\zeta}.$\\
Now for $r=\psi(x),\ s=\psi(y),\ t=\psi(z),$ the equality (3.2) becomes
$$\varphi(x)[f(y,z)]+f(x,yz)=-(\psi^*\zeta')(x,y,z)+f(x,y)+f(xy,z)$$
So
$$\overline{\psi^*\zeta'}=\overline{k}=\overline{\zeta}.$$
\end{proof}
\subsection {The equivalence between a Gr-category and a strict one}
Using Theorem 8 and Theorem of the realization of the obstruction in the group extension problem, we prove the following theorem
\begin{thm}
Each Gr-category is Gr-equivalent to a strict one.
\end{thm}
First, we prove the following lemma
\begin{lem}
Let $\mathcal{C}'$ be a strict Gr-category whose the reduced
Gr-category is $S'=(\Pi', C', \xi')$. Then, for each group
homomorphism $\psi:\Pi\to \Pi',$ there exists a strict Gr-category
$\mathcal{C},$ Gr-equivalent to the Gr-category $S=(\Pi, C, \xi),$ where $\mathcal{C}$ is regarded as a $\Pi$-module with the operation $xc=\psi(x)c, $ and $\xi$ belongs to the same cohomological class of $\psi^*{\xi'}.$
\end{lem}
\begin{proof}
We construct the strict Gr-category $\mathcal{C}$ as follows
\begin{gather*}
Ob(\mathcal{C})=\{(x, X)| \ x\in \Pi, X\in \psi(x)\} \\
Hom ((x, X), (x, Y))=\{x\}\times Hom_{\mathcal{C}'}(X,Y).
\end{gather*}
The tensor product on objects and arrows of $\mathcal{C}$ are defined by
\begin{gather*}
(x, X)\otimes (y,Y)=(xy, X\otimes Y)\\
(x,u)\otimes (y,v)=(xy, u\otimes v).
\end{gather*}
The unitivity object of $\mathcal{C}$ is $(1, I)$ where $I$ is the unitivity object of
$\mathcal{C}'.$ Readers may easily verify that the category
$\mathcal{C} $ is a strict Gr-category. Moreover, we have the isomorphisms
\[\begin{aligned}
\lambda:\Pi_0(\mathcal{C})\to \Pi\\
 \overline{(x,X)}\mapsto x
\end{aligned}\qquad\qquad
\begin{aligned}
\rho:\Pi_1(\mathcal{C})&\to \Pi_1(\mathcal{C}')=C\\
(1,c)&\mapsto c
\end{aligned}\]
and a Gr-functor $(F,\widetilde{F}):\mathcal{C} \to \mathcal{C}'$ is given by 
$$F(x,X)=X,\ \ F(x,u)=u, \ \  \Fm=id.$$
Let $\overline F=(F_0,F_1)$ be the
functor induced by $(F,\widetilde{F})$ on the reduced categories
$\mathcal S, \mathcal S',$ we have
\begin{gather*}
F_0\overline{(x,X)}=\overline{F(x,X)}=\overline X=\psi(x)\\
F_1(1,u)=\gamma_{F(1,I)}F(1,u)=\gamma_I(u)=u.
\end{gather*}
This means $F_0=\psi \lambda$ and $F_1=f,$ or $\overline F$ is the functor of the type $(\psi\lambda, f).$\\
\indent Now assume that $u$ is the associativity constraint of
$\mathcal S.$ Let $(\phi,\widetilde{\phi})$ denote the Gr-functor
from $\mathcal S$ to $\mathcal S'$ determined by the composition
$$(\phi,\widetilde{\phi} )=(G', \widetilde{G}')\circ(F,\widetilde{F})\circ(H,\widetilde {H})$$
for canonical equivalences $(H,\widetilde {H}), (G', \widetilde{G}')$. From the Proposition 2, we have $\phi = \overline F$.\\
Besides, from Theorem 6, the obstruction of the pair $(\psi\lambda, f)$ must vanish in $H^{3}(\Pi_{0}, C^{'}),$ i.e.,
$$\ \ \ (\psi\lambda)^{\ast}\xi^{'}=f_{\ast}u + \delta \widetilde{\phi}\ \ \ $$
Now if we denote $\xi=f_{\ast}u,$ the pair $J=(\lambda, f), \widetilde{J}=id$ is a Gr-functor from $\mathcal S$ to $(\Pi,C,\xi).$ Then the composition $(J,\widetilde{J})\circ(G,\widetilde{G})$ is a Gr-equivalence from $\mathcal C$ to $(\Pi,C,\xi).$
Finally, we will prove that $\xi$ belongs to the same cohomological class as $\psi^{\ast}\xi^{'}.$ Let $K=(\lambda^{-1}, f^{-1}):(\Pi, C, \xi)\rightarrow \mathcal S.$ Then $K$ together with $\widetilde{K}=id$ is a Gr-category, and the composition
$$(\phi, \widetilde{\phi})\circ (K, \widetilde{K}):(\Pi, C, \xi)\rightarrow \mathcal S'$$
is a Gr-functor. 
\[
\begin{diagram}
\node[1]{\mathcal S}
 \arrow[2]{e,t}{\Phi}
\node[2]{\mathcal S'}
\\
\node[2]{(\Pi, C, \xi)}\arrow[1]{nw,r,2}{K}
\arrow[1]{ne,r,2}{\Phi\circ K}
\end{diagram}
%\tag{1.7}
\]Clearly, $\phi\circ K$ is a functor of the type $(\psi, id)$ and therefore its obstruction vanishes, i.e., $\xi$ belongs to the same cohomological class as $\psi^{\ast}\xi^{'}.$

{\bf The proof of Theorem 9}
 
Let $\mathcal{A}$ be a Gr-category
whose the reduced Gr-category is $\mathcal{I}=(\Pi, C, \xi)$, the theorem on the realization of the obstruction (Theorem 9.2, Chapter IV [6]), there exists the group $G$ whose center is $Z(G)=C$ and such that the abstract kernel $\Pi, G, \psi)$ has the obstruction $\mathcal S'=(\Pi^{'}, C, \xi^{'}),$ where $\overline{\psi^*\xi'}=\overline{\xi}.$\\
\indent From Lemma 10, the homomorphism $\psi:\Pi\to AutG/IntG,$ defines a strict Gr-category $\C,$ Gr-equivalent to $\mathcal I=(\Pi, C, \xi).$ 
Therefore, $\C$ and $\A$ are Gr-equivalent. The theorem is completely proved.
\end{proof}
Readers may find a different proof of Theorem 9 in [10].
\begin{center}

\end{center}
Address: Department of Mathematics\\
Hanoi National University of Education\\
136 Xuan Thuy Street, Cau Giay district, Hanoi, Vietnam.\\
Email: nguyenquang272002@gmail.com
\end{document}